\RequirePackage{ifpdf}
\ifpdf 
\documentclass[pdftex]{sigma}
\else
\documentclass{sigma}
\fi

\numberwithin{equation}{section}

\begin{document}

\allowdisplaybreaks

\renewcommand{\PaperNumber}{015}

\FirstPageHeading

\renewcommand{\thefootnote}{$\star$}
\ShortArticleName{Quasi-Linear Algebras and Integrability (the
Heisenberg Picture)}

\ArticleName{Quasi-Linear Algebras and Integrability\\ (the
Heisenberg Picture)\footnote{This paper is a contribution to the
Proceedings of the Seventh International Conference ``Symmetry in
Nonlinear Mathematical Physics'' (June 24--30, 2007, Kyiv,
Ukraine). The full collection is available at
\href{http://www.emis.de/journals/SIGMA/symmetry2007.html}{http://www.emis.de/journals/SIGMA/symmetry2007.html}}}

\Author{Luc VINET~$^\dag$ and Alexei ZHEDANOV~$^\ddag$}
\AuthorNameForHeading{L.~Vinet and A.~Zhedanov}

\Address{$^\dag$~Universit\'e de Montr\'eal PO Box 6128, Station
Centre-ville, Montr\'eal QC H3C 3J7, Canada}
\EmailD{\href{mailto:luc.vinet@umontreal.ca}{luc.vinet@umontreal.ca}}

\Address{$^\ddag$~Donetsk Institute for Physics and Technology,
Donetsk 83114, Ukraine}
\EmailD{\href{mailto:zhedanov@yahoo.com}{zhedanov@yahoo.com}}

\ArticleDates{Received November 16, 2007, in f\/inal form January
19, 2008; Published online February 06, 2008}

\Abstract{We study Poisson and operator algebras with the
``quasi-linear property'' from the Heisenberg picture point of
view. This means that there exists a set of one-parameter groups
yielding an explicit expression of dynamical variables (operators)
as functions of ``time'' $t$. We show that many algebras with
nonlinear commutation relations such as the Askey--Wilson,
$q$-Dolan--Grady and others satisfy this property. This provides
one more (explicit Heisenberg evolution) interpretation of the
corresponding integrable systems.}

\Keywords{Lie algebras; Poisson algebras; nonlinear algebras;
Askey--Wilson algebra; Dolan--Grady relations}

\Classification{17B63; 17B37; 47L90}

\section{Classical version}
\setcounter{equation}{0}

Assume that there is a classical Poisson manifold with the Poisson
brackets (PB) $\{x,y\}$ def\/ined for all dynamical variables $x$,
$y$ belonging to this manifold. Of course, it is assumed that the
PB satisfy the standard conditions:
\begin{enumerate}\itemsep=0pt
\item[(i)] $\{x, \alpha_1 y+ \alpha_2 z \} = \alpha_1 \{x,y\} +
\alpha_2 \{x, z\}$ -- linearity ($\alpha_1$, $\alpha_2$ are
arbitrary constants);

\item[(ii)] $\{x,y\} = -\{y,x\}$ -- antisymmetricity;

\item[(iii)] $\{x,yz\} = z\{x,y\} + y\{x,z\}$ -- the Leibnitz
rule;

\item[(iv)] $\{\{x,y\},z\} + \{\{y,z\},x\} + \{\{z,x\},y\} =0$ --
the Jacobi identity.
\end{enumerate}

If one chooses the  ``Hamiltonian'' $H$ (i.e.\ some dynamical
variable belonging to this manifold), then we have the standard
Hamiltonian dynamics: all variables $x(t)$ become depending on
an~additional time variable $t$ and the equation of motion is
def\/ined as
\begin{gather}
 \dot x(t) = \{x(t),H\} \label{dot_x}.
 \end{gather} Of course, the
Hamiltonian $H$ is independent of $t$, because $\dot H = \{H, H\}
=0$. More generally, a~dynamical variable $Q$ is called the
integral of motion if $\dot Q =0$. Clearly, this is equivalent to
the condition $\{Q,H\}=0$.

From \eqref{dot_x} we have $\ddot x = \{\{x,H\},H\}$ and more
generally
\[
\frac{d^n x}{dt^n} = \{ \dots \{ x,H\}, H\}, \dots H\} = \{x,
H^{(n)}\},
\]
where notation $\{x, H^{(n)}\}$ means $n$ repeated PB. Using this
formula we can write down the formal Taylor expansion of an
arbitrary dynamical variable in the form
\begin{gather}
x(t) = x(0) + t \dot x(0) + \dots + \frac{t^n}{n!} \frac{d^n
x(t)}{dt^n}\Bigl|_{t=0} + \cdots  \nonumber \\
\phantom{x(t) }{} =x(0) + t \{x(0),H\} + \dots + \frac{t^n}{n!} \:
\{x(0), H^{(n)}\} + \cdots. \label{x_ser}
\end{gather}
 Assume now that there are $N$ dynamical variables
$x_1(0), x_2(0), \dots, x_N(0)$ (not depending on~$t$) with the
following property: for every $k=1,\dots,N$ the PB of the variable
$x_k(0)$ with the Hamiltonian $H$ depends {\it linearly} on all
$x_i(0)$:
\begin{gather}
 \{x_k(0),H\} = \sum_{s=1}^N F_{ks}(H) x_s(0) +
\Phi_k(H), \qquad k=1,2,\dots, N, \label{Ans_xHN}
\end{gather}
 where
$F_{ks}(H)$, $\Phi_k(H)$ are some functions of the variable $H$
only. Using rules (i)--(iv) for the PB we have
\[
\{\{x_k(0), H\}, H\} =\sum_{s=1}^N F_{ks}(H) \{x_s(0),H\} =
\sum_{s=1}^N F^{(2)}_{ks}(H) x_s(0) + \Phi^{(1)}_k(H),
\]
where
\[
F^{(2)}_{ks}(H) = \sum_{i=1}^N F_{ki}(H) F_{is}(H), \qquad
\Phi_k^{(1)}(H) = \sum_{i=1}^N F_{ki}(H) \Phi_i(H).
\]
We see that the double PB $\{\{x_k(0), H\}, H\}$ is again linear
with respect to $x_i(0)$.

In what follows we will use the matrix notation assuming that
$F(H)$ is the $N \times N$ matrix with the entries $F_{ik}(H)$ and
$\Phi(H)$ is the $N$-dimensional vector with the components
$\Phi_i(H)$. Then it is seen that $F^{(2)}(H)$ is the square of
the matrix $F(H)$ and the vector $\Phi^{(1)}$ is obtained by the
applying of the matrix $F(h)$ to the vector $\Phi(H)$:
$\Phi^{(1)}(H) = F(H) \Phi(H)$ in accordance with usual
conventions in linear algebra.

By induction, we obtain the following formula
\begin{gather}
 \{x_k(0),
H^{(n)} \} = \sum_{s=1}^N F^{(n)}_{ks}(H) x_s(0) +
\Phi^{(n-1)}_k(H), \label{x_k_Hn}
\end{gather} where the matrix $F^{(n)}(H)$
means the $n$-th power of the matrix $F_{ik}$ and the vector
$\Phi^{(n-1)}(H)$ is
\[
\Phi^{(n-1)}(H) = F^{(n-1)} \Phi(H).
\]
Now using \eqref{x_ser} we can write down the formula
\[
x_k(t) = \sum_{n=0}^{\infty} \{x_k(0), H^{(n)}\}   \frac{t^n}{n!}.
\]
 By
\eqref{x_k_Hn} we can present the above formula in the form
\begin{gather} x_k(t) = \sum_{s=1}^N E_{ks}(H;t) x_s(0) +
G_k(H;t), \label{x_k(t)}
\end{gather} where the $N \times N$ matrix $E(H;t)$ and the $N$-dimensional
vector $G(H,t)$ are def\/ined as \begin{gather} E(H;t) = \exp(t
F(H)), \qquad G(H;t) = \left( \int_{\tau=0}^t \exp(\tau F(H))d\tau
\: \right) \Phi(H). \label{EGt}
\end{gather}
We used the ordinary def\/inition of the exponential function of
the matrix:
\[
\exp(tF) = \sum_{n=0}^{\infty} \frac{t^n F^{(n)}}{n!}
\]
and
\[
\int_{\tau=0}^t \exp(\tau F)d\tau = \sum_{n=0}^{\infty}
\frac{t^{n+1} F^{(n)}}{(n+1)!}.
\]
We see from \eqref{x_k(t)} that every variable $x_k(t)$ is a
linear function with respect to initial variab\-les~$x_s(0)$.

\begin{proposition}
Assume that there is an analytical matrix function $E(H;t)$ and an
analytical vector function $G(H;t)$ such that condition
\eqref{x_k(t)} holds for every $k=1,2,\dots, N$. Then this
condition is equivalent to condition \eqref{Ans_xHN}.
\end{proposition}
Proof of this proposition is elementary.

So far, we did not concretize the choice of the Hamiltonian $H$.
Now we assume that {\it all} variables $x_s(0)$, $s=1,\dots, N$
can be chosen as Hamiltonians: if one puts $H=x_s(0)$ then we will
have linearity property \eqref{Ans_xHN} with respect to all other
variables $x_i$, $i\ne s$. It is naturally to call such algebras
the {\it quasi-linear} Poison algebras. This means that linearity
property holds for all generators $x_i$, apart from the
Hamiltonian variable~$x_s$.

Quasi-linear algebras possess a remarkable property: if one
chooses any basic variable $x_j$, $j=1,2,\dots, N$ as a
Hamiltonian then for time evolution of all other variables $x_i$,
$i=1,\dots, N$, $i\ne j$ we have linear property
\begin{gather*} x_i(t) =
{\sum_{s=1}^N}\,{\vphantom{\sum}}' \xi_{ijs}(x_j;t) x_s(0) +
\eta_{ij}(x_j;t), \qquad
i \ne j, 
\end{gather*}
where the functions $\xi_{ijs}(x_j;t)$ and $\eta_{ij}(x_j;t)$ can
be explicitly calculated as was shown above. Notation
${\sum\limits_{s=1}^N}\,{\vphantom{\sum}}'$ means that we exclude
the term with $s=j$ from the sum (note that excluding value of $s$
coincides with the f\/ixed ``number'' $j$ of the Hamiltonian).

In this case we obtain strong restrictions for the expression of
the Poisson brackets $\{x_i,x_k\}$. Indeed, from this property we
have that for any pair of variables $x_i$, $x_k$ one has the
relation
\begin{gather}
 \{x_i,x_k\} = {\sum_{s=1}^N}\,{\vphantom{\sum}}' F_{iks}(x_k) x_s +
\Phi_{ik}(x_k), \qquad i,k =1,2, \dots, N, \qquad i \ne k
\label{Pois_ik}
\end{gather} with some functions $F_{iks}(x_k)$ and
$\Phi_{ik}(x_k)$. These functions cannot be taken arbitrarily.
Indeed, we can also write down the similar condition but with
reversed order of variables $x_i$, $x_k$:
\begin{gather*}
 \{x_k,x_i\} =
{\sum_{s=1}^N}\,{\vphantom{\sum}}' F_{kis}(x_i) x_s +
\Phi_{ki}(x_i), \qquad i,k =1,2,
\dots, N, \qquad i \ne k. 
\end{gather*} But
$\{x_k,x_i\}=-\{x_i,x_k\}$, so we have the system of conditions
\begin{gather}
 {\sum_{s=1}^N}\,{\vphantom{\sum}}' F_{kis}(x_i) x_s + \Phi_{ki}(x_i) +
{\sum_{s=1}^N}\,{\vphantom{\sum}}' F_{iks}(x_k) x_s +
\Phi_{ik}(x_k)=0 \label{Pois_cond}
\end{gather} which should be valid for {\it any} pair $x_i$,
$x_k$, $i\ne k$.

Moreover, the Jacobi identity \begin{gather}
  \{\{x_i,x_k\},x_j\} +
\{\{x_k,x_j\},x_i\} + \{\{x_j,x_i\},x_k\} =0 \label{Jac_ikj}
\end{gather} should be valid for any triple $x_i$, $x_k$, $x_j$
with distinct values $i$, $j$, $k$. In our case the Jacobi
identity \eqref{Jac_ikj} is reduced to a system of
functional-dif\/ferential equations \begin{gather} W_{ijk} +
W_{jki} + W_{kij}=0, \label{Jac_W} \end{gather} where
\[
W_{ijk} = \{\{x_i,x_k\},x_j\} = \Phi'_{ik}(x_k) \{x_k,x_j\} +
{\sum_{s=1}^N}' \left(F_{iks}(x_k) \{x_s,x_j\} + F'_{iks}(x_k) x_s
\{x_k, x_j\}  \right),
\]
where $F'(x)$ means the derivative of the function $F(x)$.

Conditions \eqref{Pois_cond} and \eqref{Jac_W} can be considered
as a system of nontrivial functional-dif\/ferential equations for
unknown functions $F_{iks}(x_k)$ and $\Phi_{ik}(x_k)$.

Of course, there is a trivial solution of these conditions when
all functions are constants, i.e.\ $F_{iks}(x_k)$ and
$\Phi_{ik}(x_k)$ do not depend on their arguments~$x_k$. If,
additionally, $\Phi_{ik} \equiv 0$ then we obtain well-known
Lie--Poisson algebras with commutation relations
\[
\{x_i, x_k\} = \sum_{s=1}^N c_{iks} x_s
\]
with the structure constants $c_{iks}$ satisfying standard
restrictions following from \eqref{Pois_cond} and \eqref{Jac_W}
(see, e.g.~\cite{KarMa}).

Consider the most general f\/inite-dimensional Poisson algebras
with nonlinear Poison brackets \begin{gather} \{x_i,x_k\}=
h_{ik}(x_1, \dots, x_N), \label{nlin_PB} \end{gather} where
$h_{ik}(x_1, \dots, x_N)$ are smooth functions of $N$ variables
$x_1, \dots, x_N$. For general theory connected with these
algebras see e.g.~\cite{KarMa}. In order for the variables $x_i$,
$x_k$ to form a Poisson algebra, the functions~$h_{ik}$ should
satisfy some strong restrictions. In particular, in~\cite{FG}
necessary and suf\/f\/icient conditions were obtained in case when
$h_{ik}(x_1, \dots, x_N)$ are quadratic functions. We will provide
several examples when Poisson algebras of type~\eqref{nlin_PB}
satisfy conditions~\eqref{Pois_cond} and~\eqref{Jac_W} and hence
are quasi-linear. General classif\/ication of all Poisson algebras
with such property is an interesting open problem.

Note that in~\cite{KarMa} the so-called semi-linear Poisson
algebras (as well as their operator analogues) were considered.
Such algebras are in general nonlinear but they have linearity
property with respect to some prescribed generators. However the
semi-linear algebras introduced in~\cite{KarMa} do not possess, in
general, the property~\eqref{Pois_ik} and hence they are not
quasi-linear algebras in our sense.

\section[''Quantum'' (operator) version]{``Quantum'' (operator) version}

\setcounter{equation}{0} Now assume that $X_k$, $k=1,2,\dots$ are
operators which act on some linear space (either f\/inite or
inf\/inite-dimensional).

The time dynamics is def\/ined by the Heisenberg equations:
\begin{gather*} \dot
B =[H,B], 
\end{gather*} where $H$ is an operator called the
Hamiltonian and $[H,B] \equiv HB-BH$ is the commutator.

As in the previous section, assume that there exist $N$ operators
$X_k$ and a Hamiltonian $H$ such that the conditions
\begin{gather}
 [H,X_k] =
\sum_{s=1}^N F_{ks}(H) X_s + \Phi_k(H), \label{X_k_H_cond}
\end{gather} where
$F_{ks}(H)$, $\Phi_k(H)$ are some functions depending only on the
Hamiltonian $H$.

Introduce the $n$-th repeated commutator
\[
[H,[H,\dots, [H,B]\dots]={\rm ad}_H^n B.
\]
As in the previous section we have
\begin{gather*}
{\rm  ad}_H^n X_k=\sum_{s=1}^N
F^{(n)}_{ks}(H) X_s + \Phi^{(n-1)}_k(H), 
\end{gather*} where
the matrix $F^{(n)}(H)$ means the $n$-th power of the matrix
$F_{ik}$ and the vector $\Phi^{(n-1)}(H)$ is
\[
\Phi^{(n-1)}(H) = F^{(n-1)} \Phi(H).
\]
The Heisenberg time evolution is described by the one-parametric
group in a standard manner \begin{gather} X_k(t) = \exp(tH) X_k
\exp(-tH)
 = X_k + t [H,X_k] + \dots + \frac{t^n}{n!}\, {\rm
ad}_H^n X_k + \cdots. \label{X_k(t)} \end{gather}
 Note that in our approach the
time variable $t$ can be an arbitrary complex parameter and the
``Hamiltonian'' $H$ need not be a Hermitian operator. In physical
applications it is usually assumed that $H$ is a Hermitian
operator and we need to change $t \to it$ in order to obtain the
usual Heisenberg picture where the time~$t$ is a real parameter.

Again, as in classical case we obtain that $X_k(t)$ is a linear
combination of initial operators~$X_s$: \begin{gather}
 X_k(t) = \sum_{s=1}^N
E_{ks}(H;t) X_s + G_k(H;t), \label{op_X_k(t)} \end{gather} where
the $N \times N$ operator-valued matrix $E(H;t)$ and the
$N$-dimensional operator-valued vector $G(H,t)$ are def\/ined by
the same formulas~\eqref{EGt}.

Note that if the operators $X_s$ form a Lie algebra:
\[
[X_i,X_k]=\sum_{s=1}^N g_{ik}^s X_s,
\]
where $g_{ik}^s$ are the structure constants, then any operator
$X_j=H$ taken as a Hamiltonian, satisf\/ies conditions
\eqref{X_k_H_cond} with $F_{ks}$ the constants (not depending on
$H$) and $\Phi_k=0$, hence we have formula \eqref{X_k(t)} with
$G_k =0$ which is a standard action of the one-parameter Lie group
corresponding to the generator $H=X_j$:
\begin{gather*} \exp(tX_j) X_k
\exp(-tX_j) =\sum_{s=1}^N E_{ks}(t) X_s. 
\end{gather*} Of
course, in the case of the Lie algebra we can construct full Lie
group due to linearity property ($g_{ik}^s$ are constants).
However, if $F_{ks}(H)$ depend on an operator $H$ then we see that
only a set of one-parameter groups exists with the linearity
property~\eqref{X_k(t)}. Our next problem therefore, will be how
to construct algebras satisfying the property~\eqref{X_k_H_cond}
for dif\/ferent possible choices of the operator~$H$.

The simplest possibility is the same as in the previous section:
we demand that {\it any} opera\-tor~$X_j$ can be taken as a
Hamiltonian $H=X_j$, $j=1,2,\dots,N$. We call such algebras the
quasi-linear operator algebras. Condition~\eqref{X_k_H_cond} is
replaced with a system of conditions \begin{gather} [X_j,X_k] =
{\sum_{s=1}^N}\,{\vphantom{\sum}}' F_{jks}(X_j) X_s +
\Phi_{jk}(X_j), \qquad k,j= 1,2 \dots, N, \qquad k \ne j.
\label{X_k_j_cond} \end{gather} Then for all possible choices of
the Hamiltonian $H$ the time evolution (the Heisenberg picture)
has the linearity property with respect to all operators apart
from the Hamiltonian: \begin{gather*} X_k^{(j)}(t)=\exp(tX_j) X_k
\exp(-tX_j) ={\sum_{s=1}^N}\,{\vphantom{\sum}}' E_{jks}(X_j;t) X_s
+ G_{jk}(X_j;t),
\end{gather*} where the functions $E_{jks}(X_j;t)$ and
$G_{jk}(X_j;t)$ are easily calculated from the functions
$F_{jks}(X_j)$ and $\Phi_{jk}(X_j)$.

Of course, compatibility analysis of conditions \eqref{X_k_j_cond}
as well as checking of the Jacobi identity
\[
[[X_i,X_k],X_j] + [[X_k,X_j],X_i] + [[X_j,X_i],X_k] =0, \qquad
i,j,k =1,2, \dots, N
\]
is rather a nontrivial problem even for the case when all
operators $X_k$ are f\/inite-dimensional. Instead, we present
concrete examples of the operator algebras possessing quasi-linear
property.

However f\/irst we need to generalize our scheme to include the
so-called ``extension'' operators.

\section{Quasi-linear algebras with extension}

Assume that we have a set of $N$ so-called ``basic'' operators
$X_1, X_2, \dots X_N$ and for any opera\-tor~$X_i$, $i=1,2,\dots,
N$ we have also ``a tower'' of extension, i.e.\  a set of $M_i \ge
0$ operators $Y^{(i)}_k$, $k=1,2,\dots, M_i$, such that we have
the commutation relations \begin{gather} [X_i,X_k] =
{\sum_{s=1}^N}\,{\vphantom{\sum}}' F_{iks}(X_i) X_s +
\sum_{s=1}^{M_i}G_{iks}(X_i) Y_s^{(i)} + \Phi_{ik}(X_i), \nonumber\\
i,k=1,2,\dots, N , \qquad i\ne k \label{ext_XX} \end{gather} and
\begin{gather} [X_i,Y_k^{(i)}] =
{\sum_{s=1}^N}\,{\vphantom{\sum}}' U_{iks}(X_i) X_s +
\sum_{s=1}^{M_i}V_{iks}(X_i) Y_s^{(i)} + W_{ik}(X_i), \nonumber
\\ i=1,2,\dots, N, \qquad k=1,2,\dots, M_i \label{ext_XY} \end{gather} with some
functions $F_{iks}(x), \dots, W_{ik}(x)$. In what follows we will
assume that all these functions are polynomials.

Thus the extension operators $Y_k^{(i)}$ enter the linear
combinations in commutators \eqref{ext_XX}, \eqref{ext_XY} but
only operators $X_i$ can be chosen as ``Hamiltonians''. In
general, commutation relations between operators $Y_k^{(i)}$  are
not def\/ined as well as  the commutation relations between the
Hamiltonian $X_i$ and operators $Y_k^{(j)}$ from another ``tower''
(i.e.\ when $j \ne i$).

From \eqref{ext_XX}, \eqref{ext_XY} it follows immediately that
repeated commutators with the ``Hamilto\-nians''~$X_i$ have the
similar structure
\[
{\rm ad}_{X_i}^n X_k = {\sum_{s=1}^N}\,{\vphantom{\sum}}'
F_{iks}^{(n)}(X_i) X_s + \sum_{s=1}^{M_i}G_{iks}^{(n)}(X_i)
Y_s^{(i)} + \Phi_{ik}^{(n)}(X_i), \qquad i,k=1,2,\dots, N
\]
and
\[
{\rm ad}_{X_i}^n Y_k^{(i)} ={\sum_{s=1}^N}\,{\vphantom{\sum}}'
U_{iks}^{(n)}(X_i) X_s + \sum_{s=1}^{M_i}V_{iks}^{(n)}(X_i)
Y_s^{(i)} + W_{ik}^{(n)}(X_i)
\]
with polynomials $F_{iks}^{(n)}(x), \dots, W_{ik}^{(n)}(x)$ which
can be obtained from $F_{iks}(x), \dots, W_{ik}(x)$ by obvious
explicit procedures.

Hence we have an explicit Heisenberg evolution picture:
\[
e^{tX_i} X_k e^{-tX_i} ={\sum_{s=1}^N}\,{\vphantom{\sum}}' \tilde
F_{iks}(X_i;t) X_s + \sum_{s=1}^{M_i} \tilde G_{iks}(X_i;t)
Y_s^{(i)} + \tilde \Phi_{ik}(X_i;t), \qquad i,k=1,2,\dots, N
\]
and
\[
e^{tX_i} Y_k^{(i)} e^{-tX_i}={\sum_{s=1}^N}\,{\vphantom{\sum}}'
\tilde U_{iks}(X_i;t) X_s + \sum_{s=1}^{M_i} \tilde V_{iks}(X_i;t)
Y_s^{(i)} + \tilde W_{ik}(X_i;t),
\]
where the functions $\tilde F_{iks}(x;t), \dots, \tilde
W_{ik}(x;t)$ are obtained from the polynomials $F_{iks}^{(n)}(x),
\dots$, $W_{ik}^{(n)}(x)$ by an obvious way, e.g.
\[
\tilde F_{iks}(x;t) = \sum_{n=0}^{\infty} \frac{t^n}{n!} \,
F_{iks}^{(n)}(x)
\]
etc.

Thus, as in the previous section, we have an explicit {\it linear}
Heisenberg evolution of the operators $X_i$, $Y_k^{(i)}$ under the
action of~$N$ Hamiltonians $X_i$, $i=1,2,\dots,N$. In contrast to
the previous section, the operators $Y_k^{(i)}$ do not in general
provide the linear evolution, i.e.\ they cannot be chosen as
Hamiltonians. They serve only as an auxiliary tool in our picture.

The total number $M$ of extension operators $Y_k^{(i)}$ can be
less then $M_1+M_2 + \dots + M_N$, because in special cases some
of these operators can coincide, as we will see later.

It will be convenient to call the corresponding algebras the
quasi-linear algebras  of the type $(N,M)$, where~$N$ is number of
the ``true Hamiltonians'' $X_1, \dots, X_N$ and~$M$ is a total
number of extension operators~$Y_k^{(i)}$. For quasi-linear
algebras without extension we will use the symbol~$(N,0)$.

Of course, essentially the same picture is valid in classical case
if one replaces the commutators $[\cdots]$ with the Poisson
brackets $\{\cdots\}$.

In the next sections we construct simplest examples of the
quasi-linear algebras in both classical and quantum (operator)
pictures.

\section[The $q$-oscillator algebra]{The $\boldsymbol{q}$-oscillator algebra}

 Consider the simplest case when we have
only two operators $X$, $Y$ (i.e.\ $N=2$) and we demand that property
\eqref{X_k_H_cond} will be satisf\/ied when $H$ is chosen as
either $X$ or $Y$, i.e.\ in this case we have the quasi-linear
algebra of type $(2,0)$ without extension. From previous
considerations it follows that two conditions \begin{gather}
[X,Y]= F_1(Y) X + \Phi_1(Y), \qquad [Y,X]=F_2(X) Y + \Phi_2(X)
\label{op_2-case} \end{gather} should be valid, where $F_i(z)$,
$\Phi_i(z)$, $i=1,2$ are some functions.

In the classical case we will have similar conditions:
\begin{gather} \{x,y\}= F_1(y) x + \Phi_1(y), \qquad
\{y,x\}=F_2(x) y + \Phi_2(x) \label{cl_2-case} \end{gather} with
commutators $[\cdots]$ replaced with the Poisson brackets
$\{\cdots\}$. We see that necessary condition for compatibility
of~\eqref{op_2-case} or~\eqref{cl_2-case} is \begin{gather} F_1(Y)
X + \Phi_1(Y) + F_2(X) Y + \Phi_2(X)=0. \label{2-cond}
\end{gather} In the classical case (i.e.\ when $X$, $Y$ are two
independent commuting variables) we have a~simple functional
equation for 4 unknown functions $F_i(z)$, $\Phi_i(z)$, $i=1,2$.
Assume f\/irst that the function~$F_2(x)$ is not a constant. Then
taking two distinct values $x=\alpha$ and $x=\beta$ and
considering \eqref{2-cond} as a system for unknowns $F_1(y)$,
$\Phi_1(y)$ we easily f\/ind that these functions should be linear
in~$y$. Quite similarly, we f\/ind that functions $F_2(x)$,
$\Phi_2(x)$ should be linear as well. After simple calculations we
obtain in this case the most general solution \begin{gather}
\{x,y\} = \alpha xy + \beta_1 x + \beta_2 y + \gamma,
\label{cl_qosc} \end{gather} where $\alpha$, $\beta_1$, $\beta_2$,
$\gamma$ are arbitrary constants ($\alpha \ne 0$). By an
appropriate af\/f\/ine transformation of the variables $x \to
\xi_1 x + \eta_1$, $ y \to \xi_2 y + \eta_2$ we can reduce
\eqref{cl_qosc} to a canonical form
\begin{gather*} \{x,y\}=\alpha xy -1 
 \end{gather*} in the case if the
polynomial $\alpha xy + \beta_1 x + \beta_2 y + \gamma$ is
irreducible. If this polynomial is reducible (i.e.\ can be
presented as a product of two linear polynomials in~$x$,~$y$) then
we
have the canonical form \begin{gather*} \{x,y\}= \alpha xy. 
\end{gather*}
In quantum case we have correspondingly either algebra of the form
\begin{gather} [X,Y] =\alpha XY + \gamma, \label{q_osc1} \end{gather} if $\gamma \ne 0$,
or the Weyl operator pair \begin{gather*} XY = qYX . 
\end{gather*} Note that
relation~\eqref{q_osc1} can be rewritten in  other canonical form
\begin{gather} XY - qYX =1 \label{q-osc} \end{gather} with some
parameter $q \ne 0, 1$. The pair of operators $X$, $Y$ satisfying
the commutation relation~\eqref{q-osc} is called the
$q$-oscillator algebra.

In this case the Heisenberg evolution becomes very simple. Indeed,
consider, e.g.~$Y$ as the Hamiltonian. We have from \eqref{q-osc}
\[
[Y,[Y,X]]=\omega^2 Y^2 X - \omega Y,
\]
where $\omega = 1-q$, and, by induction,
\[
ad_Y^n X = \omega^n Y^n X - \omega ^{n-1} Y^{n-1}, \qquad
n=1,2,\dots.
\]
Hence we have explicitly
\[
e^{tY} X e^{-tY} = e^{\omega t Y} X - \phi(Y;t),
\]
where
\[
\phi(y;t)= \int_0^t e^{\omega \tau y} d \tau = \frac{e^{\omega
ty}-1}{\omega y}.
\]
When $\omega \to 0$ (i.e.\ $q \to 1$) one obtains
\[
e^{tY} X e^{-tY} = X-t + \omega t Y(X-t/2) + O(\omega^2).
\]
Here the f\/irst term $X-t$ in the right-hand side corresponds to
the simple Weyl shift of the harmonic oscillator. Indeed, when
$X$, $Y$ satisfy the Heisenberg--Weyl commutation relation $[X,Y]=1$
then
\[
e^{tY} X e^{-yY}= X-t,
\]
i.e.\ the operator $X$ is shifted by a constant under action of
the Hamiltonian $H=Y$. The second term $\omega t Y(X-t/2)$
describes a small ``$q$-deformation'' of the Weyl shift.

\section[The Askey-Wilson algebra]{The Askey--Wilson algebra}

 The Askey--Wilson algebra $AW(3)$
\cite{GLZ,ZAW} can be presented in several equivalent forms. One
of them consists of 3 operators $K_1$, $K_2$, $K_3$ with the
commutation relations
\begin{gather}  [K_1, K_2]=K_3, \label{AW_1} \\
 [K_2, K_3] = 2\rho K_2 K_1 K_2 + a_1 \{K_1, K_2\} + a_2 K_2^2 +
c_1 K_1 + d K_2 + g_1, \nonumber \\  [K_3, K_1] = 2\rho K_1 K_2
K_1 + a_2 \{K_1, K_2\} + a_1 K_1^2 + c_2 K_2 + d K_1 + g_2
\nonumber \end{gather} with some constants $\rho$, $a_1$, $a_2$,
$c_1$, $c_2$, $d$, $g_1$, $g_2$. The symbol $\{\cdots\}$ stands
for the anticommutator: $\{X,Y\} \equiv XY+YX$.

The $AW(3)$ algebra belongs to quasi-linear type of $(2,1)$.
Indeed, we can rewrite it in the form
\begin{gather} [K_1, K_2]=K_3,
\label{AW_11} \\  [K_2, [K_2,K_1]] = R_2(K_2) K_1 + R_1(K_2) K_3 +
R_0(K_2), \nonumber
\\  [K_1, [K_1, K_2]] = S_2(K_1) K_2 + S_1(K_1) K_3 +
S_0(K_1), \nonumber
\end{gather} where
\begin{gather*}
R_2(x) = -2 \rho x^2 - 2 a_1 x - c_1, \qquad R_1(x) = -2 \rho x -
a_1, \qquad R_0(x) = -a_2 x^2 -d x - g_1,
\\
S_2(x) = -2 \rho x^2 - 2 a_2 x - c_2, \qquad S_1(x) = 2 \rho x +
a_2, \qquad S_0(x) = -a_1 x^2 -d x - g_2.
\end{gather*}
We see that $AW(3)$ algebra in the form \eqref{AW_11} has
quasi-linear structure with $N=2$, $M=1$. The operators $K_1$,
$K_2$ play the role of the true Hamiltonians whereas the operator
$K_3$ is the only extension.

If $\rho \ne 0$ one can present the same algebra in an equivalent
(more symmetric) form in terms of operators $X$, $Y$, $Z$
\cite{ZAW, Ter,WiZa} \begin{gather}
 XY - q YX = Z + C_3, \qquad YZ - q ZY = X
+ C_1, \qquad ZX - q XZ = Y + C_2, \label{q_AW}
\end{gather}
 or, in a pure
commutator form \begin{gather*} [X,Y] = (q-1) YX + Z +C_3, \qquad
[Y,Z] = (q-1)
ZY+ X +C_1, \nonumber \\  [Z,X] = (q-1) XZ + Y +C_2. 
\end{gather*} In this form we have quasi-linear algebra of type $(3,0)$ {\it
without extension}, i.e.\ in this case all 3~operators $X$, $Y$,
$Z$ of the algebra play the role of the Hamiltonians. Hence in
this case the linearity property holds for all operators $X$, $Y$,
$Z$.

It is instructive to derive directly transformation property of,
say, operator $X$ under Hamiltonian action of the operator $Y$. We
have
\begin{gather}
{\rm ad}_Y^n X = U_n(Y)X + V_n(Y) Z + W_n(Y), \label{ad_n_AW}
\end{gather}
where $U_n(x)$, $V_n(x)$, $W_n(x)$ are polynomials in $x$. Initial
conditions are obvious:
\[
U_0(x) =1, \qquad V_0(x) = W_0(x)=0.
\]
Derive recurrence relations for the polynomials $U_n(x)$,
$V_n(x)$, $W_n(x)$. For this goal let as apply the operator ${\rm
ad}_Y$ with respect to \eqref{ad_n_AW}:
\begin{gather*}
{\rm ad}_Y^{n+1} X = U_n(Y)[Y,X] + V_n(Y) [Y,Z] =
\left((1-q)YU_n(Y) + q^{-1} V_n(Y)\right)X \\
\phantom{{\rm ad}_Y^{n+1} X =}{} + \left((1-q^{-1})Y V_n(Y)
-U_n(Y) \right)Z -C_3 U_n(Y) + q^{-1} C_1 V_n(Y).
\end{gather*}
Hence we have
\begin{gather}
 U_{n+1} = (1-q)YU_n + q^{-1} V_n, \qquad V_{n+1} = -U_n + (1-q^{-1})YV_n, \nonumber \\
 W_{n+1}=-C_3 U_n + C_1 q^{-1} V_n. \label{rec_UVW}
 \end{gather}
From recurrence relations \eqref{rec_UVW} it is possible to f\/ind
polynomials $U_n$, $V_n$, $W_n$ explicitly. In particular, it is
clear that $U_n(Y)$ is a polynomial of degree $n$ and $V_n(Y)$,
$W_n(Y)$ are polynomials of degree $n-1$.

For the transformed operator $X(t)$ we have
\[
X(t) = e^{tY}Xe^{-tY} = E_1(t;Y) X + E_2(t;Y) Z + E_0(t;Y),
\]
where
\[
E_1(t;Y) = \sum_{n=0}^{\infty} \frac{t^n U_n(Y)}{n!}, \qquad
E_2(t;Y) = \sum_{n=0}^{\infty} \frac{t^n V_n(Y)}{n!}, \qquad
E_0(t;Y) = \sum_{n=0}^{\infty} \frac{t^n W_n(Y)}{n!}.
\]
From \eqref{rec_UVW} we have a system of linear dif\/ferential
equations for unknown functions $E_1(t)$, $E_2(t)$
\[
\dot E_1(t;x) =  (1-q) x E_1 + q^{-1} E_2, \qquad \dot E_2(t;x) =
- E_1 + (1-q^{-1}) x E_2.
\]
Initial conditions are $E_1(0)=1$, $E_2(0)=0$. This system has
constant coef\/f\/icients (not depending on $t$) and hence can be
elementary integrated:
\[
E_1(t;x) = a_{11}(x) e^{\omega_1(x)t} + a_{12}(x)
e^{\omega_2(x)t}, \qquad E_2(t;x) = a_{21}(x) e^{\omega_1(x)t} +
a_{22}(x) e^{\omega_2(x)t},
\]
where $\omega_{1,2}(x)$ are roots of the characteristic equation
\[
\omega^2 + x(2-q-q^{-1})(x-\omega) + q^{-1}=0
\]
and the coef\/f\/icients $a_{ik}(x)$ can be found by standard
methods.

The function $E_0(t;x)$ can then be found as
\[
E_0(t;x) = -C_3 \: \int_0^{t} E_1(\tau;x) d\tau + C_1 q^{-1} \:
\int_0^{t} E_2(\tau;x) d\tau.
\]
We see that the Heisenberg evolution is described by elementary
functions (linear combinations of exponents) in the argument~$t$.
But the coef\/f\/icients in these linear combinations depend on
operator~$H$.

Similar expressions are valid for all other possible choices of
the Hamiltonian (i.e.\ $H=X$ or $H=Z$) due to symmetric form of
algebra~\eqref{q_AW}.

Thus for $\rho \ne 0$ essentially the same $AW(3)$ algebra can be
presented in two equivalent forms~-- either as an algebra of type
$(2,1)$ with the only extension or as an algebra of type $(3,0)$
without extension. However if $\rho=0$ then the only type $(2,1)$
is known (this special case corresponds to the so-called quadratic
Racah algebra QR(3)~\cite{GLZ}).

So far, we have no general results in quantum (i.e.~operator) case
concerning classif\/ication scheme. In particular, we do not know,
whether $AW(3)$ algebra is the only quasi-linear algebra of type
$(2,1)$ or $(3,0)$. In the classical case the situation is
slightly better: we already showed that for the type $(2,0)$ the
only $q$-oscillator algebra appears. In the next section we
consider a~classif\/ication scheme for the type $(3,0)$.

\section{Classical case of type (3,0)}

 Consider the classical case of
quasi-linear algebras of type $(3,0)$. This means that we have
3~dynamical variables $x$, $y$, $z$ satisfying the Poisson bracket
relations \begin{gather*}  \{y, z\} = F_1(x,y,z), \qquad \{z, x\}
= F_2(x,y,z),
\qquad \{x, y\} = F_3(x,y,z), 
\end{gather*} where the functions
$F_i(x,y,z)$ should be chosen in such way to satisfy the Jacobi
identity \begin{gather*}   \{\{x,y\},z\} + \{\{y,z\},x\} +
\{\{z,x\},y\} =0.
\end{gather*} It is easily seen that this condition is
equivalent to the relation~\cite{FG} \begin{gather} ({\bf F},
\mbox{rot}\, {\bf F})=0, \label{F_rot} \end{gather} where the
vector $\bf F$ has Cartesian components $(F_1, F_2, F_3)$ and
$\mbox{rot}$ is standard dif\/ferentiation curl operator acting on
the vector~$\bf F$.

Relation \eqref{F_rot} has an obvious solution
\begin{gather} {\bf F} = \nabla
Q(x,y,z), \label{F_nabla} \end{gather} where $\nabla$ is the
gradient operator and $Q(x,y,z)$ is a function of 3~variables. In
this case we have the def\/ining Poisson relations in the form
\begin{gather} \{y, z\} = Q_x, \qquad \{z, x\} = Q_y, \qquad \{x,
y\} = Q_z, \label{Pois_Q} \end{gather} where $Q_x$ means
derivation with respect to~$x$ etc. Note that in this case the
function $Q(x,y,z)$ is the Casimir element of the algebra, i.e.\
$\{x,Q\} = \{y,Q\}=\{z,Q\}=0$. Poisson brackets of such types are
sometimes called the Nambu brackets or Mukai--Sklyanin
algebras~\cite{OdeRu}. However, the Nambu bracket \eqref{F_nabla}
do not exhaust all admissible Poisson algebras with 3~generators.
There are non-trivial examples corresponding to $\mbox{rot}\, {\bf
F} \ne 0$.

Consider restrictions on the functions $F_i(x,y,z)$ coming from
the quasi-linear property.

On the one side, choosing $y$ to be a Hamiltonian, we have by this
property \begin{gather} \{x,y\} = F_3(x,y,z) = \Phi_1(y) x +
\Phi_2(y) z + \Phi_3(y). \label{xy_1} \end{gather}

On the other side, choosing $x$ to be a Hamiltonian, we have
analogously \begin{gather}
 \{x,y\} = \Phi_4(x) y + \Phi_5(x) z + \Phi_6(x)
\label{xy_2} \end{gather} with some functions $\Phi_i(x)$,
$i=1,2,\dots,6$.

From \eqref{xy_1} and \eqref{xy_2} we see that $F_3(x,y,z)$ should
be a polynomial having degree no more than 1 with respect to each
variable $x$, $y$, $z$ and the most general form of this
polynomial is
\[
 F_3(x,y,z) = \alpha_3 xy + \beta_{31} x + \beta_{32} y +
\beta_{33}z + \gamma_3.
\] Quite analogously, taking functions
$F_1$, $F_2$ we obtain \begin{gather} F_i(x,y,z) = \alpha_i x_k
x_l + \sum_{s=1}^3 \beta_{is}x_s + \gamma_i, \label{F_i}
\end{gather} where $x_1=x$, $x_2=y$, $x_3 =z$ and $\alpha_i$,
$\beta_{is}$, $\gamma_i$ are some constants. As usual in rhs of
\eqref{F_i} notation~$x_k x_l$ means that triple $(i,k,l)$ has no
coinciding entries.

Condition \eqref{F_rot} imposes strong restrictions upon the
coef\/f\/icients $\alpha_i$, $\beta_{ik}$, $\gamma_i$. It is
convenient to analyze the canonical forms of obtained algebras up
to af\/f\/ine transformations $x_i \to \xi_i x_i + \eta_i$ with
some constants $\xi_i$, $\eta_i$ (of course we demand that $\xi_1
\xi_2 \xi_3 \ne 0$).

(i) If diagonal entries $\beta_{ii}$ are all  nonzero then
necessarily $\alpha_1 =\alpha_2 = \alpha_3 =\alpha$. We will
assume $\alpha \ne 0$ (otherwise we will obtain the Lie--Poisson
algebras). Then the matrix $\beta$ should be symmetric:
$\beta_{ik}=\beta_{ki}$. These conditions are also suf\/f\/icient
for validity of the Jacobi identity. After an appropriate
af\/f\/ine transformation we reduce our algebra to the canonical
form \begin{gather*}
F_i(x,y,z) = \alpha x_k x_l + x_i + \gamma_i 
\end{gather*}
with only 4 independent free parameters: $\alpha$ and $\gamma_i$,
$i=1,2,3$. This is exactly the classical version of the
Askey--Wilson algebra $AW(3)$ \cite{KoZhe1}.

(ii) If one of diagonal entries is zero, say $\beta_{33}=0$, then
again we have the same necessary and suf\/f\/icient conditions
$\alpha_1 =\alpha_2 = \alpha_3 =\alpha \ne 0$ and
$\beta_{ik}=\beta_{ki}$. The canonical form is \begin{gather*} F_1
= \alpha x_2 x_3 + x_1 + \gamma_1, \qquad F_2 = \alpha x_1 x_3 +
x_2 +
\gamma_2, \qquad F_3 = \alpha x_2 x_1 + \gamma_3. 
\end{gather*} In
both cases (i) and (ii) we have the Nambu--Poisson  brackets
\eqref{Pois_Q} with \begin{gather*} Q= \alpha xyz + (x^2 + y^2 +
z^2)/2 +
\gamma_1 x + \gamma_2 y + \gamma_3 z 
\end{gather*} for the
case $(i)$ and \begin{gather*} Q= \alpha xyz + (x^2 + y^2)/2 +
\gamma_1 x +
\gamma_2 y + \gamma_3 z 
\end{gather*} for the case $(ii)$.

(iii) Two diagonal entries are zero, say $\beta_{22}=
\beta_{33}=0$ and $\beta_{11} \ne 0$. Then there exist
2~possibilities: in the f\/irst case  $\alpha_1 =\alpha_2 =
\alpha_3 =\alpha \ne 0$ and the canonical form of the algebra is
\begin{gather*} F_1 = \alpha x_2 x_3 + x_1 + \gamma_1, \qquad F_2
= \alpha x_1 x_3 +
\gamma_2 , \qquad F_3 = \alpha x_2 x_1 +\gamma_3. 
\end{gather*}
This is again the algebra of Nambu type with \begin{gather*} Q=
\alpha xyz +
x^2/2 + \gamma_1 x + \gamma_2 y + \gamma_3 z. 
\end{gather*}

(iv) As in the previous case we have $\beta_{22}= \beta_{33}=0$
and $\beta_{11} \ne 0$. But now there is a~possibility of
noncoinciding entries $\alpha_i$. Namely we have $\alpha_2=
\alpha_3$ but $\alpha_1 \ne \alpha_2$. We will assume that
$\alpha_2 \alpha_1 \ne 0$. Then we have the canonical form
\begin{gather*} F_1 = \alpha_1 x_2 x_3 + x_1, \qquad F_2 =
\alpha_2 x_1 x_3 , \qquad F_3 =
\alpha_2 x_2 x_1.  
\end{gather*}

This case does not belong to the Nambu type.

(v) There is further degeneration of the previous case. Namely we
can allow $\alpha_1 =0$, $\alpha_2 \ne 0$ or $\alpha_2 =0$,
$\alpha_1 \ne 0$. In the f\/irst case we have the canonical form
\begin{gather*} F_1 = x_1 + \gamma_1, \qquad F_2 = \alpha_2 x_1
x_3 , \qquad F_3 = \alpha_2 x_2 x_1.
\end{gather*} In the second case the canonical form is \begin{gather*} F_1
= \alpha_1 x_2 x_3 + x_1 + \gamma_1, \qquad F_2 = x_3 , \qquad F_3
= x_2.
\end{gather*}

(vi) Assume that all diagonal entries are zero
$\beta_{11}=\beta_{22} = \beta_{33}=0$. Then if $\alpha_1 \alpha_2
\alpha_3 \ne 0$ we have the canonical form \begin{gather*} F_1 =
\alpha_1 x_2x_3, \qquad F_2 = \alpha_2 x_1 x_3 , \qquad F_3 =
\alpha_3 x_2 x_1.
\end{gather*}

It is interesting to note that obtained algebras of types
(i)--(vi) correspond to the list of so-called 3-dimensional
skew-polynomial algebras introduced by Bell and
Smith~\cite{Smith,Rosenberg}. In our case we have classical
(i.e.~Poisson brackets) analogues of the 3-dimensional
skew-polynomial algebras.

Concerning other aspects and applications of the classical
(Poissonic) version of the $AW(3)$ algebra see e.g.~\cite{KoZhe1}.

\section{Beyond the AW-algebra}

\setcounter{equation}{0} Consider the so-called Dolan--Grady
relations \cite{DG} for two operators $A_0, A_1$ \begin{gather}
[A_0,[A_0,[A_0,A_1]]] = \omega^2   [A_0,A_1], \qquad
[A_1,[A_1,[A_1,A_0]]] = \omega^2   [A_1,A_0], \label{DG}
\end{gather} where $\omega$ is an arbitrary constant. These
relations generate the so-called inf\/inite-dimensional Onsager
algebra~\cite{Davies} which plays a crucial role in algebraic
solution of the Ising model as well as of some more general models
in statistical physics~\cite{Perk}.

On the other hand, we can consider the DG-relations from the
``quasi-linear'' algebras point of view.  Indeed, introduce the
operators \begin{gather} A_2=[A_0,A_1], \qquad A_3=[A_0,
[A_0,A_1]], \qquad A_4 = [A_1, [A_1,A_0]]. \label{DG_ext}
\end{gather} Then we see that
\begin{gather*}
{\rm ad}_{A_0}^{2n+2} A_1 = \omega^{2n} A_3, \qquad {\rm
ad}_{A_0}^{2n+1} A_1 = \omega^{2n} A_2, \qquad n=0,1,2,\dots,
\\
{\rm ad}_{A_1}^{2n+2} A_0 = \omega^{2n} A_4, \qquad {\rm
ad}_{A_1}^{2n+1} A_0 = -\omega^{2n} A_2, \qquad n=0,1,2,\dots,
\\
{\rm ad}_{A_0}^{2n+2} A_2 = \omega^{2n+2} A_2, \qquad {\rm
ad}_{A_0}^{2n+1} A_2 = \omega^{2n} A_3, \qquad n=0,1,2,\dots,
\\
{\rm ad}_{A_1}^{2n+2} A_2 = \omega^{2n+2} A_2, \qquad {\rm
ad}_{A_1}^{2n+1} A_2 = -\omega^{2n} A_4, \qquad n=0,1,2,\dots.
\end{gather*}

From these formulas we obtain an explicit Heisenberg evolution of
the operators $A_0$, $A_1$, $A_2$ if the operators $A_0$ or $A_1$
are chosen as Hamiltonians:
\begin{gather*}
e^{tA_0} A_1 e^{-tA_0}= A_1 + \frac{\sinh(\omega t)}{\omega} A_2
+\frac{\coth( \omega t)-1}{\omega^2} A_3,
\\
e^{tA_0} A_2 e^{-tA_0}= \cosh(\omega t) A_2 + \frac{\sinh(\omega
t)}{\omega} A_3,
\\
e^{tA_1} A_0 e^{-tA_1}= A_0 - \frac{\sinh(\omega t)}{\omega} A_2
+\frac{\coth(t \omega)-1}{\omega^2} A_4,
\\
e^{tA_1} A_2 e^{-tA_1}= \cosh(\omega t) A_2 - \frac{\sinh(\omega
t)}{\omega} A_4.
\end{gather*}
We see that the operators $A_0$, $A_1$ can be chosen as
Hamiltonians whereas the operators $A_2$, $A_3$, $A_4$ are the
extension. Hence the DG-relations provide an example of a
quasi-linear algebra of the type $(2,3)$.

Def\/ine now the operator
\[
W= \alpha A_0 + \beta A_1 + \gamma A_2
\]
with arbitrary parameters $\alpha$, $\beta$, $\gamma$. It is seen
that
\begin{gather*}
e^{tA_1} W e^{-tA_1} = \alpha A_0 + \beta A_1 + \left(\gamma
\cosh(\omega t) - \alpha \frac{\sinh(\omega t)}{\omega} \right)
A_2 \\
\phantom{e^{tA_1} W e^{-tA_1} =}{} + \left(\alpha
\frac{\cosh(\omega t) -1}{\omega^2}   - \gamma \frac{\sinh(\omega
t)}{\omega} \right) A_4.
\end{gather*}
Analogously
\begin{gather*}
e^{\tau A_0} W e^{-\tau A_0} = \alpha A_0 + \beta A_1 +
\left(\gamma \cosh(\omega \tau) + \beta \frac{\sinh(\omega
\tau)}{\omega} \right) A_2\\
\phantom{e^{\tau A_0} W e^{-\tau A_0} =}{}  + \left(\beta
\frac{\cosh(\omega \tau) -1}{\omega^2}   + \gamma
\frac{\sinh(\omega \tau)}{\omega} \right) A_3,
\end{gather*}
where $t$ and $\tau$ are arbitrary parameters. If one chooses
\[
\alpha = \omega \gamma  \coth(\omega t/2), \qquad \beta = -\omega
\gamma  \coth(\omega \tau/2),
\]
then the terms containing $A_3$, $A_4$ disappear and we have
\[
e^{tA_1} W e^{-tA_1} = \alpha A_0 + \beta A_1 - \gamma A_2 =
e^{\tau A_0} W e^{-\tau A_0},
\]
whence
\[
TW T^{-1} =W,
\] where
\[
T = e^{-\tau A_0} e^{t A_1}.
\]
Equivalently, this means that the operator $W$ commutes with the
operator $T$:
\[
TW =WT.
\]
As was noted by Davies~\cite{Davies} this commutation relation was
crucial in Onsager's solution of the Ising model, where the
operators $T$ play the role of the transfer matrix. We see that
this relation follows directly from the quasi-linear property of
the DG-relations.

There is an obvious generalization of the DG-relations preserving
the quasi-linear property.

Indeed, def\/ine again basic operators $A_0$, $A_1$ and their
extensions $A_2$, $A_3$, $A_4$ by \eqref{DG_ext}. But now the
triple commutators ${\rm ad}_{A_0}^3 A_1$ and ${\rm ad}_{A_1}^3
A_0$ can contain not only operator $A_2$ as in DG-case~\eqref{DG}
but arbitrary linear combinations of the type
\begin{gather} [A_0,[A_0,[A_0,A_1]]] = g_1(A_0)A_1 +
g_2(A_0)A_2 + g_3(A_0)A_3 + g_0(A_0) \label{gen_DG1}
\end{gather} and
\begin{gather}
[A_1,[A_1,[A_1,A_0]]] = f_1(A_1)A_0 + f_2(A_1)A_2 + f_3(A_1)A_4 +
f_0(A_1) \label{gen_DG2}
\end{gather} with some polynomials $g_i(x)$, $f_i(x)$,
$i=0,\dots,3$. We will call these relations the generalized
DG-relations. The ordinary DG-relations corresponds to the choice
$g_0=f_0=g_1=f_1=g_3=f_3=0$ and $g_2=-f_2=\omega^2$.

From \eqref{gen_DG1} and \eqref{gen_DG2} it is obvious that for
any positive integer $n$ the repeated commutator ${\rm ad}_{A_0}^n
A_1$ has the same structure
\[
{\rm ad}_{A_0}^n A_1 = g_1^{(n)}(A_0)A_1 + g_2^{(n)}(A_0)A_2 +
g_3^{(n)}(A_0)A_3 + g_0^{(n)}(A_0)
\]
and similarly
\[
{\rm ad}_{A_1}^n A_0 = f_1^{(n)}(A_1)A_0 + f_2^{(n)}(A_1)A_2 +
f_3^{(n)}(A_1)A_4 + f_0^{(n)}(A_1)
\]
with polynomials $g_i^{(n)}(x)$ and $f_i^{(n)}(x)$ can be easily
obtained from $g_i(x)$ and $f_i(x)$. We thus have the ``Heisenberg
solvability'' property
\[
\exp(A_0 t)A_1 \exp(-A_0 t) = G_1(A_0;t)A_1 + G_2(A_0;t)A_2 +
G_3(A_0;t)A_3 + G_0(A_0;t),
\]
where functions $G_i(x;t)$ can be easily calculated in the same
manner as in~\eqref{op_X_k(t)}. (Similar relation holds for the
Heisenberg evolution of the operator~$A_1$.)

Consider a special example of these generalized DG-relations. The
so-called ``tridiagonal algebra'' proposed by Terwilliger is
generated by the two relations \cite{ITT} \begin{gather} [A_0,
A_0^2A_1 +A_1 A_0^2 - \beta A_0A_1A_0 - \gamma (A_0A_1 + A_1 A_0)
- \alpha A_1]=0 \label{Ter1}
\end{gather} and \begin{gather} [A_1, A_1^2A_0 +A_0 A_1^2 - \beta
A_1A_0A_1 - \gamma_1 (A_0A_1 + A_1 A_0) - \alpha_1 A_0]=0,
\label{Ter2}
\end{gather} where $\beta$, $\gamma$, $\gamma_1$, $\alpha$, $\alpha_1$
are some constants.

The tridiagonal algebra is closely related with the Askey--Wilson
algebra $AW(3)$. Indeed, let the operators $A_0$, $A_1$ belong to the
AW(3) algebra of the type \eqref{AW_1}. Identify $A_0 = K_1$,
$A_1=K_2$, $A_2 = K_3$ and rewrite relations \eqref{AW_1} in an
equivalent form getting rid of the operator $K_3=A_2$:
\begin{gather} A_1^2A_0 + A_0 A_1^2 +
2(\rho -1) A_1 A_0 A_1 + a_1(A_1A_0+ A_0 A_1)  \nonumber \\
\qquad{}+a_2 A_1^2 + c_1 A_0 + d A_1 + g_1 =0 \label{1AW}
\end{gather} and
\begin{gather}
A_0^2A_1 + A_1 A_0^2 + 2(\rho -1) A_0 A_1 A_0 + a_2(A_1A_0+ A_0
A_1) \nonumber
\\ \qquad{}+a_1 A_0^2 + c_2 A_1 + d A_0 + g_2 =0. \label{2AW}
\end{gather}
 Now
applying the operator ${\rm ad}_{A_0}$ to  \eqref{2AW} we obtain
relation \eqref{Ter1}. Analogously, applying the operator ${\rm
ad}_{A_1}$ to \eqref{1AW} we obtain relation \eqref{Ter2}.
Parameters of the tridiagonal algebra~\eqref{Ter1},~\eqref{Ter2}
are related with parameters of the $AW(3)$ algebra as follows
\[
\beta= 2(1-\rho), \qquad \alpha = -c_2, \qquad \alpha_1 = -c_1,
\qquad \gamma=-a_2, \qquad \gamma_1= -a_1.
\]
Thus the tridiagonal algebra \eqref{Ter1}, \eqref{Ter2} follows
from the $AW(3)$ algebra~\eqref{AW_1}. However the reciprocal
statement is not valid: the tridiagonal algebra is larger than
$AW(3)$~\cite{ITT}.

It is easily verif\/ied that the tridiagonal algebra \eqref{Ter1},
\eqref{Ter2} can be presented in an equivalent form:
\begin{gather}
 [A_0,[A_0,[A_0,A_1]]] = (2-\beta)(A_0 A_3- A_0^2 A_2) + 2
\gamma A_0 A_2 - \gamma A_3 + \alpha A_2, \nonumber \\
 [A_1,[A_1,[A_1,A_0]]] = (2-\beta)(A_1 A_4 + A_1^2 A_2) -
2\gamma_1 A_1 A_2 - \gamma_1 A_4 - \alpha_1 A_2,
\label{semiT}\end{gather} where the operators $A_2$, $A_3$, $A_4$
are def\/ined as \eqref{DG_ext}. Comparing~\eqref{semiT}
with~\eqref{gen_DG1} and~\eqref{gen_DG2} we see that the
tridiagonal algebra~\eqref{Ter1} indeed belongs to the
quasi-linear type with the coef\/f\/icients
\begin{gather*}
g_0=g_1=0, \qquad g_2(x)=(\beta-2)x^2 + 2\gamma x + \alpha, \qquad
g_3(x) = (2-\beta) x - \gamma,
\\
f_0=f_1=0, \qquad f_2(x) = (2-\beta) x^2 - 2 \gamma_1 x -
\alpha_1, \qquad f_3(x) = (2-\beta)x - \gamma_1.
\end{gather*}
As in the case of the ordinary DG-relations, the tridiagonal
algebra is the quasi-linear algebra with two basic operators
$A_0$, $A_1$ and 3~extensions $A_2$, $A_3$, $A_4$, i.e.\ it has
the same type $(2,3)$ as the DG-algebra.

In a special case $\gamma=\gamma_1=0$ we obtain the so-called
$q$-deformation of the Dolan--Grady relations \cite{Ter,Ter1}:
\begin{gather}
 [A_0,[A_0,[A_0,A_1]_q]_{q^{-1}}] = \alpha
[A_0,A_1], \qquad [A_1,[A_1,[A_1,A_0]_q]_{q^{-1}}] = \alpha_1
[A_1,A_0], \label{qDG} \end{gather} where $[X,Y]_q = q^{1/2}XY -
q^{-1/2}YX$ is so-called $q$-commutator. The parameter $q$ is
connected  with the parameter $\beta$ by the relation $\beta= q+
q^{-1}.$ This deformed DG-algebra plays an important  role in
theory of quantum $XXZ$ Heisenberg model, Azbel--Hofstadter model
etc \cite{Ba,BasKo}.

It would be interesting to study what is the meaning of the exact
solvability of the Heisenberg picture for the $q$-deformed
DG-relations \eqref{qDG} in corresponding exactly solvable models.

\section{Quasi-linear algebras and exactly solvable systems}

Recently Odake and Sasaki proposed an interesting approach to
exactly solvable classical and quantum mechanical models in the
Heisenberg picture \cite{OS, OS1}. They noticed that almost all
exactly solvable one-dimensional quantum models admit exact
solution in both Schr\"odinger and Heisenberg pictures.

Consider, e.g. the one-dimensional model described by the standard
one-dimensional Hamiltonian \begin{gather} H = p^2/2 + U(x),
\label{1_Ham} \end{gather} where $p$ is the momentum operator and
$U(x)$ is a potential. Of course, the standard Heisenberg
commutation relation between coordinate and momentum operators
$[x,p]=i$ is assumed.

Introduce the operator $X= f(x)$ with some function $f(x)$ and
consider the Heisenberg evolution of this operator under the
Hamiltonian $H$: \begin{gather*} X(t) = \exp(iHt) X \exp(-iHt) = X
+ it [H,X] - \frac{t^2}{2}   [H,[H,X]]+   \dots +
\frac{(it)^n}{n!} \,{\rm ad}_H^n X + \cdots. 
\end{gather*}
In general for arbitrary Hamiltonian $H$ the structure of the
expression ${\rm ad}_H^n X$ will be too complicated, and so it is
hopeless to f\/ind explicit solution in the Heisenberg picture.
Assume, however, that the Hamiltonian $H$ and the function $f(x)$
satisfy the restriction
\begin{gather}
 [H,[H,X]]={\rm ad}_H^2 X = g_1(H) X + g_2(H) Y + g_0(H),
\label{SO_Ans} \end{gather}
 with some functions $g_i(H)$, $i=0,1,2$, where
we introduce the operator
\[
Y=[H,X] = -2i f'(x) - f''(x).
\]  Then
it is obvious that ansatz~\eqref{SO_Ans} (proposed by Odake and
Sasaki \cite{OS,OS1}) leads to an explicit solution
\begin{gather*} X(t)
= G_1(H;t) X + G_2(H;t) Y + G_0(H;t), 
\end{gather*}
where the functions $G_i(H;t)$, $i=0,1,2$ have elementary behavior
in time $t$ (i.e.\ they can be expressed in terms of trigonometric
or hyperbolic functions).

Simple considerations (see \cite{OS}) lead to the conclusion that
$g_1(H)$ and $g_0(H)$ should be linear polynomials in $H$ whereas
$g_2(H)$ should be a constant such that $2g_2 = -
g_1'(H)=-\alpha_1$ with some constant $\alpha_1$.

It is easy to f\/ind (see  \cite{OS}) that the most general
function $f(x)$ satisfying this property is a~solution of the
equation
\[
f''(x) = -\tfrac{1}{2}(\alpha_1 f(x) + \beta_1)
\]
with arbitrary constants $\alpha_1$, $\beta_1$. Hence $f(x)$ is an
elementary function: it is either a quadratic polynomial in $x$,
or a superposition of two exponents (plus a constant).
Corresponding potentials~$U(x)$ can also be easily
found~\cite{OS}~-- they coincide with well-known
``exactly-solvable'' potentials in quantum mechanics:

(i) singular oscillator:
\[
U(x) = a_1 x^2 + a_2 x^{-2} + a_0;
\]

(ii) the Morse potential
\[
U(x) = a_1 e^{-2x} + a_2 e^{-x} + a_0;
\]

(iii) the P\"oschl--Teller potential
\[
U(x) = a_1 \sin^{-2}x + a_2 \cos^{-2} x+ a_0
\]
(in the latter case one can replace trigonometric functions with
hyperbolic ones that gives hyperbolic analogs of the
P\"oschl--Teller potential).

Observe now that two operators $H$, $X$ together with the third
operator $Y=[H,X]$ form some quasi-linear algebra. In this case
the operator $Y$ plays the role of extension. Indeed, we already
established the f\/irst commutation relation \eqref{SO_Ans} which
in our case can be rewritten in the form \begin{gather}
 [H,Y]= (\alpha_1 H +
\alpha_0)X - \frac{\alpha_1}{2} Y + \beta_1 H + \beta_0.
\label{Jac_1} \end{gather} The second commutation relation is
verif\/ied easily: \begin{gather} [Y,X]= - f'^2(x) =
\frac{\alpha_1}{2} X^2 + \beta_1 X + \varepsilon \label{Jac_2}
\end{gather} with some constant $\varepsilon$. Relations
\eqref{Jac_1} and \eqref{Jac_2} show that we have a quasi-linear
algebra of type $(2,1)$ with two Hamiltonians $H$, $X$ and the
only extension~$Y$. This algebra is equivalent to so-called
quadratic Jacobi algebra $QJ(3)$ which was considered
in~\cite{GLZ} as a hidden symmetry algebra of exactly solvable
Schr\"odinger Hamiltonians (see also~\cite{LV}). The Jacobi
algebra is a~special case of the algebra $AW(3)$~\cite{GLZ}. It
corresponds to the following choice of the parameters
$\rho=a_2=d=0$ in \eqref{AW_1}. We thus see that the Odake--Sasaki
approach \cite{OS,OS1} for the Schr\"odinger
Hamiltonians~$H$~\eqref{1_Ham} leads naturally to the Jacobi
algebra~$QJ(3)$.

Consider a dif\/ference analogue of the Schr\"odinger Hamiltonian.
We choose \begin{gather} H = A(s) T^+ + C(s) T^- + B(s),
\label{AW_H} \end{gather} where operator $H$ acts on the space of
functions $f(s)$ depending on a variable $s$ and
opera\-tors~$T^{\pm}$ are the standard shift operators
\[
T^{\pm} f(s) = f(s \pm 1).
\]
The operator $H$ is a second-order dif\/ference operator. Hence
the
eigenvalue equation \begin{gather*} H f(s) = \lambda f(s) 
\end{gather*} can
be considered as a dif\/ference analogue of the corresponding
one-dimensional Schr\"odinger equation.

Introduce also the operators $X$ which acts as a multiplication
\[
X f(s) = x(s) f(s)
\]
with some function $x(s)$ which will be called a ``grid'' and the
commutator $Y = [H,X]$.

We now would like to f\/ind when the relation \begin{gather} {\rm
ad}_H^2 X = [H,Y] =W_1(H) X + W_2(H)Y + W_0(H) \label{AW_Ans}
\end{gather} holds, where $W_i(H)$, $i=0,1,2$ are some polynomials
in~$H$. If \eqref{AW_Ans} is fulf\/illed then obviously we have an
explicit time dynamics of the ``grid'' operator $X$ under the
Heisenberg evolution with the Hamiltonian $H$:
\[
X(t) \equiv \exp(Ht) X \exp(-Ht) = Q_1(H;t) X + Q_2(H;t) Y +
Q_0(H;t)
\]
with functions $Q_i(H;t)$ which can be easily explicitly
calculated from~$W_i(H)$.

We will assume that $A(s) C(s) \ne 0$ and the grid $x(s)$ is
nondegenerated, i.e.~$x(s) \ne x(s+1)$ and $x(s) \ne x(s+2)$.
Under such conditions it is always possible to assume that $C(s)
=1$. Indeed, using similarity transformation $H \to F(s) H
F^{-1}(s)$ with some functions~$F(s)$ one can achieve the
condition $C(s) =1$. Such similarity transformation does not
change the operator relation~\eqref{AW_Ans}.

Assume that polynomial $W_2(H)$ have a degree $N$ with respect to
the variable $H$. Then it is clear from nondegeneracy of $x(s)$
that degrees of polynomials $W_1(H)$ and $W_0(H)$ cannot
exceed~$N+1$. Assuming that ${\rm deg}\,(W_1) = {\rm
deg}\,(W_0)=N+1$ (otherwise we again obtain degeneracy of $x(s)$)
we obtain from \eqref{AW_Ans} an operator identity of the form
\[
\sum_{k=-N-1}^{N+1} E_k(s)T^k =0
\]
with coef\/f\/icients $E_k(s)$ which can be explicitly calculated
after substitution of \eqref{AW_H} and $X=x(s)$
into~\eqref{AW_Ans}. From this identity we should have
\[
E_k(s) \equiv 0
\]
for all $k=-N-1, -N, \dots, N+1$. Assume that $N \ge 2$. Then the
highest-order conditions $E_{-N-1}(s)=E_{N+1}(s) \equiv 0$ are
reduced to very simple f\/irst-degree recurrence relations for the
grid $x(s)$:
\[
\xi_Nx(s+1)  + \eta_N x(s) + \zeta_N =0
\]
and
\[
\xi_Nx(s-1)  + \eta_N x(s) + \zeta_N =0,
\]
where $\xi_N$, $\eta_N$, $\zeta_N$ are some complex
coef\/f\/icients not depending on~$s$. It is easily verif\/ied
that for all possible choices of the coef\/f\/icients these
conditions are incompatible with nondegeneracy conditions $x(s)
\ne x(s+1)$ and $x(s-1) \ne x(s+1)$ for the grid~$x(s)$. Hence
necessarily $N \le 1$. Thus $W_2(H)$ should be a linear polynomial
(or a constant) in $H$ and $W_0(H)$, $W_1(H)$ should be quadratic
(or linear) polynomials in $H$. Consider again the highest-order
conditions $E_2(s) = E_{-2}(s) \equiv 0$:
\[
\xi x(s+1) + \eta x(s) + x(s-1) + \zeta =0
\]
and
\[
\xi x(s-1) + \eta x(s) + x(s+1) + \zeta =0
\]
with some constants $\xi$, $\eta$, $\zeta$.

These conditions are compatible with non-degeneracy conditions for
the grid $x(s)$ only if $\xi=1$ and then we obtain the linear
recurrence relation \begin{gather*} x(s+1) + x(s-1) + \eta x(s) +
\zeta =0.
\end{gather*}
 This equation is well known -- it describes the
so-called ``Askey--Wilson grids'' (AW-grid)~\cite{AW, NSU}. The
most general is the ``$q$-quadratic'' grid:
\[
x(s) = c_1 q^s + c_2 q^{-s} + c_0
\]
with some constants $c_0$, $c_1$, $c_2$ and a ``basic'' parameter
$q$. In the limit $q=1$ we obtain the quadratic grid:
\[
x(s) = c_2 s^2 + c_1 s + c_0.
\]
Further simple analysis of the relation~\eqref{AW_Ans} allows to
obtain explicit expressions for~$A(s)$ and~$B(s)$. It appears (we
omit technical details) that corresponding operator $H$ coincides
with the dif\/ference Askey--Wilson operator on the grid $x(s)$
(for details see \cite{AW,NSU,VZB}). Thus we showed that the only
condition \eqref{AW_Ans} concerning integrability in the
Heisenberg picture of the second-order dif\/ference operator leads
uniquely to the Askey--Wilson dif\/ference operator and
corresponding AW-grid. It can be compared with results
of~\cite{VZB}, where the similar statement was obtained but under
dif\/ferent conditions.

It is then directly verif\/ied that the second (``dual'')
algebraic relation \begin{gather} [Y,X] = V_1(X)H + V_2(X) Y +
V_0(X) \label{dual_AW} \end{gather} holds automatically with
$\deg(V_2(X)) \le 1$ and $\deg(V_{0,1}(X)) \le 2$.

Combining \eqref{AW_Ans} and \eqref{dual_AW} we immediately
conclude that operators $H$, $X$ together with their extension
operator $Y=[H,X]$ form the (generic) Askey--Wilson algebra
$AW(3)$ \eqref{AW_11}. Thus the $AW(3)$ algebra can be derived
uniquely from the ``Heisenberg solvability'' principle. This
results seems to be new. Note that in \cite{OS} and \cite{OS1} it
was directly verif\/ied that the Askey--Wilson Hamiltonian
\eqref{AW_H} (as well as all its special cases corresponding to
classical orthogonal polynomials) provides an exact solution for
the ``grid variable'' $X$ in the Heisenberg picture. In our
approach these results follow directly from the quasi-linear
property of the $AW(3)$ algebra.

Concerning classical (i.e.~Poisson brackets) analogue of the $AW(3)$
algebra and corresponding exactly solvable systems see
e.g.~\cite{ZheKo}.

\section{Conclusion}

 We demonstrated that the ``quasi-linear''
algebras (with possible extension) can be considered as a natural
generalization of the Lie algebras: they preserve the linearity
property with respect to one-parameter subgroups (exponential
mapping) constructed from the f\/ixed set of the ``Hamiltonians''.
This provides an exact time evolution (Heisenberg picture) with
respect to all these Hamiltonians.  The situation is almost the
same in both classical (with Poisson brackets instead of
commutators) and quantum picture. Many exactly solvable models in
classical and quantum mechanics admit an alternative description
in terms of a corresponding quasi-linear algebra.

There are many open questions and problems connected with
quasi-linear algebras:

(i) is it possible to give complete classif\/ication of all
f\/inite-dimensional quasi-linear algebras (with or without
extension)? We were able to construct such classif\/ication only
in the classical case (i.e.\ Poisson algebras) for the simplest
cases of types $(2,0)$ and $(3,0)$ without extensions. We hope
that in quantum (operator) case such classif\/ication is possible
at least in the case when all operators are f\/inite-dimensional.

(ii) in our def\/inition of the quasi-linear algebras, e.g.\ in
condition~\eqref{X_k_H_cond} the polynomials $F$
 (depending on the Hamiltonians $H$) stand to the left side of the operators $X_k$.
 Of course it is possible to def\/ine the ``right'' quasi-linear algebras with
 the property
\begin{gather} [H,X_k] = {\sum_{s=1}^N} X_s F_{ks}(H)  + \Phi_k(H),
\label{X_k_H_cond_r} \end{gather} instead of \eqref{X_k_H_cond}.
It is easily verif\/ied that for all considered examples of the
quasi-linear algebras the ``right'' version coincides with the
``left'' version (with possible modif\/ication of the structure
parameters). Is this property valid in general situation? We hope
that the answer is positive.

Note that if the operators $H$ and all $X_i$ are Hermitian and if
$F_{ks}(x)$ and $\Phi_k(x)$ are polynomials in $x$ (as happens
e.g.\ for already considered examples of integrable systems in
quantum mechanics) then obviously the ``left'' and the ``right''
versions are equivalent; indeed, taking Hermitian conjugation of
\eqref{X_k_H_cond_r} we have
\begin{gather*}
[H,X_k]^{\dagger} = -[H,X_k] = {\sum_{s=1}^N}  F_{ks}^*(H) X_s  +
\Phi_k^*(H),
\end{gather*}
($F^*$ means complex conjugation of corresponding polynomials),
i.e.\ we obtained the ``left'' version from the initial ``right''
version. However, from general mathematical point of view we
cannot assume that operators $H$ and $X_i$ are Hermitian, so the
problem of equivalence between the ``left'' and the ``right''
version remains open.

(iii) possible applications to exactly solvable models. As we
know, the Jacobi algebra $QJ(3)$ describes exactly solvable
one-dimensional quantum Hamiltonians with the singular oscillator,
Morse and P\"oschl--Teller potentials. The Askey--Wilson algebra
$AW(3)$ describes all ``classical'' second-order dif\/ference
equations. By ``classical'' we mean that these equations have
solutions in terms of classical orthogonal polynomials on
non-uniform grids (the most general are the Askey--Wilson
dif\/ference equations~\cite{VZB}). What about multi-dimensional
generalization of these results? An interesting approach was
proposed in~\cite{Ba}, where the $q$-deformed DG-relations are
applied to integrable models like $XXZ$-chain on the one side and
to some generalization of the Askey--Wilson polynomials on the
other side.

(iv) the ``quasi-linear'' algebras do not describe all interesting
non-linear algebraic objects in theory of integrable systems.
Among such algebras are so-called Sklyanin algebras introduced
in~\cite{Skl} and then generalized in many papers. The Sklyanin
algebra consists from 4~generators $S_i$, $i=0,\dots,4$ with some
special quadratic commutation relations between all possible pairs
of the generators. These algebras are closely related with
elliptic solutions of the Yang--Baxter equation for the 8-vertex
model in statistical physics. Recently it was shown that the
Sklyanin algebras play an important role in theory of biorthogonal
rational functions on elliptic grids~\cite{Ros1,SpZhS}. Already on
the classical level (Poisson bracket version) the time evolution
with respect to any ``Hamiltonian'' $S_i$ is described by elliptic
functions \cite{KoZhe}. Hence in this case time evolution does not
satisfy linearity property with respect to generators. It would be
very interesting to f\/ind a~basic property of these algebras
which generalize the quasi-linear property.

\subsection*{Acknowledgments}

A.Zh.\ thanks Centre de Recherches Math\'ematiques of the
Universit\'e de Montr\'eal for hospitality and T.~Ito, A.~Kiselev,
M.~Nesterenko and P.~Terwilliger for discussions. The authors
would like to thank referees for valuable remarks and comments.

\pdfbookmark[1]{References}{ref}
\LastPageEnding

\end{document}